\documentclass[letterpaper, 10 pt, conference]{IEEEtran}%
\usepackage{graphics}
\usepackage{epsfig}
\usepackage{cite}
\usepackage{amsmath}
\usepackage{amssymb}
\usepackage{mathrsfs}
\usepackage{amsfonts}
\usepackage{graphicx}%
\providecommand{\U}[1]{\protect\rule{.1in}{.1in}}
\IEEEoverridecommandlockouts
\newtheorem{innercustomthm}{Theorem}
\newenvironment{customthm}[1]
  {\renewcommand\theinnercustomthm{#1}\innercustomthm}
  {\endinnercustomthm}

\newtheorem{innercustomdef}{Definition}
\newenvironment{customdef}[1]
  {\renewcommand\theinnercustomdef{#1}\innercustomdef}
  {\endinnercustomdef}  

\newtheorem{innercustomass}{Assumption}
\newenvironment{customass}[1]
  {\renewcommand\theinnercustomass{#1}\innercustomass}
  {\endinnercustomass}  
  
\newtheorem{innercustomcor}{Corollary}
\newenvironment{customcor}[1]
  {\renewcommand\theinnercustomcor{#1}\innercustomcor}
  {\endinnercustomcor}  
  
\newtheorem{innercustomrem}{Remark}
\newenvironment{customrem}[1]
  {\renewcommand\theinnercustomrem{#1}\innercustomrem}
  {\endinnercustomcor}   

\begin{document}

\title{
Optimal Control for Speed Harmonization of Automated Vehicles
}

\author{Andreas A. Malikopoulos, {\itshape{Senior Member, IEEE}}, Seongah Hong, B. Brian Park  {\itshape{Member, IEEE}}, \\Joyoung Lee, and Seunghan Ryu 
\thanks{Andreas A. Malikopoulos is with the Department of Mechanical Engineering, University of Delaware, DE 19716 USA (phone: 302-831-2889; e-mail: 
        {andreas@udel.edu)}.}
\thanks{Seongah Hong, B. Brian Park, and Seunghan Ryu are with the Department of Civil and Environmental Engineering, University of Virginia Charlottesville, VA 22904-4742 USA USA (phone: 434-924-6347; e-mail: 
        {sh3zm@virginia.edu; bp6v@virginia.edu; sr5ae@virginia.edu)}.}       
\thanks{Joyoung Lee is with the Department of Civil Engineering, New Jersey Institute of Technology, Newark, NJ 07102 (phone: 973-596-2475; e-mail: 
        {jo.y.lee@njit.edu)}.}
 }


\maketitle 
\begin{abstract}
This article addresses the problem of controlling the speed of a number
of automated vehicles before they enter a speed reduction
zone on a freeway. 
We formulate the control problem and provide an analytical, closed-form solution that can be implemented in real time. The
solution yields the optimal acceleration/deceleration of each vehicle
under the hard safety constraint of rear-end collision avoidance. The effectiveness of the solution is evaluated through a microscopic simulation testbed
and it is shown that the proposed approach significantly reduces 
both fuel consumption and travel time. In particular,  for three different traffic volume
levels, fuel consumption for each vehicle is reduced by 19-22\% compared to the baseline scenario, in which human-driven vehicles are considered, by 12-17\% compared to the variable speed limit algorithm, and by 18-34\% compared to the vehicular-based speed harmonization (SPD-HARM) algorithm. Similarly, travel time is improved by 26-30\% compared to the baseline scenario, by 3-19\% compared to the VSL algorithm, and by 31-39\% compared to the vehicular-based SPD-HARM algorithm.
 
\end{abstract}
\thispagestyle{empty} \pagestyle{empty}

\indent


\section{Introduction}

\subsection{Motivation}

In a rapidly urbanizing world, we need to make fundamental transformations
in how we use and access transportation. We are currently witnessing an increasing integration of our energy, transportation, and cyber networks, which, coupled with the human interactions, is giving rise to a new level of complexity in transportation \cite{Malikopoulos2015}. 
As we move to increasingly complex transportation systems \cite{Malikopoulos2016c},
new control approaches are needed to optimize the impact on system
behavior of the interaction between vehicles at different traffic scenarios. Intersections, merging roadways, speed reduction zones along with
the drivers' responses to various disturbances 
are the primary sources of bottlenecks that contribute to traffic
congestion \cite{Malikopoulos2016a}. In 2015, congestion caused
people in urban areas in US to spend 6.9 billion hours more on the road
and to purchase extra 3.1 billion gallons of fuel, resulting in
a total cost estimated at \$160 billion \cite{Schrank2015}.

Speed harmonization (SPD-HARM) is one of the major applications operated in the US towards reducing congestion. SPD-HARM encompasses a group of strategies that are intended to form a series of coherent traffic stream along the roadway by regulating traffic speeds, especially when there is a need of reducing speed down the road. The fundamental idea of such SPD-HARM is to mitigate the loss of highway performance by preventing traffic break-downs and keeping bottlenecks operating at constant traffic feeds. Such idea is well represented in the fundamentals of traffic flow theory: a traffic break-down at the bottleneck can be prevented by progressively guiding the upstream traffic to equal the downstream traffic flow, so the upstream traffic smoothly runs into the downstream traffic and can pass through the bottleneck without disruptions. Using such insight, numerous SPD-HARM strategies were developed with various approaches to determine the speed policy for the system. The control algorithms for SPD-HARM may include traditional intelligent transportation systems technologies or use information provided to connected and automated vehicles (CAVs) to enforce speed on individual vehicle - the individual vehicle-based SPD-HARM can be viewed as an extreme form of VSL. It is important to note that the vehicular-based SPD-HARM should be differentiated from any type of cruise control system such as adaptive cruise control (ACC) or cooperative ACC (CACC), since the latter are intended to facilitate the inter-vehicular interactions by maintaining a given distance between the vehicles as opposed to a strategy for vehicles to approach a bottleneck. 

In this paper, we develop a framework that allows each vehicle to optimally control its speed before entering a bottleneck. The objective is to derive the optimal acceleration/deceleration of each vehicle within a ``control zone" of appropriate length right before entering a speed reduction zone. The latter causes bottleneck that builds up as vehicles exceed the bottleneck capacity. By optimizing the vehicles' acceleration/deceleration in the upstream, the time of reaching the speed reduction zone is controlled optimally, and thus the recovery time from the congested area is minimized. In addition, the vehicle can minimize stop-and-go driving, thereby conserving momentum and energy. Eliminating the vehicles' stop-and-go driving aims at minimizing transient engine operation, and thus we can have direct benefits in fuel consumption \cite{Malikopoulos2008b}. 

The contributions of this paper are: (1) the problem formulation to control optimally the speed of a number of automated vehicles cruising on a freeway before they enter a speed reduction zone, (2) the analytical, closed-form solution of this problem along with a rigorous analysis that reveals the conditions under which the rear-end collision avoidance constraint does not become active within the control and speed reduction zones, and (3) bringing the gap between the existing algorithms and the real-world implications by providing an optimal solution that can be implemented in a vehicle in real time. We should emphasize that our approach focuses on controlling optimally the vehicle before its entry to the speed reduction zone and not on inter-vehicle interactions, and threfore, it is fundamentally different from the ACC and CACC approaches.

\subsection{Literature Review}

Traditionally, the SPD-HARM has been realized through variable message signs (VMS), variable speed limit (VSL) and the rolling SPD-HARM (a.k.a., pace-car technique) \cite{roberts_i-70_2012}. Both
VMS and VSL systems employ the display gantries mounted along roadways
to deliver messages or control schemes. Another method for SPD-HARM is 
the rolling SPD-HARM, which uses designated patrol vehicles
entering the traffic to hold a traffic stream at a lower speed, and thus, traverse the congestion area smoothly while mitigating
shock waves.

The application of SPD-HARM
has been mainly evolved through VSL which appeared
to be more effective and efficient than VMS
and the rolling SPD-HARM \cite{robinson_examples_2000,roberts_i-70_2012}.
State-of-the-art VSL method employs a proactive approach applying a control action beforehand and then anticipate the behavior
of the system (traffic) \cite{khondaker_variable_2015}. Even
though this proactive approach has made VSL a popular method over the years, it provides a sub-optimal solution since it is based on a heuristic approach \cite{frejo_global_2012}.

SPD-HARM methods can be broadly categorized into the (1) reactive approach and (2) proactive approach. The reactive approach initiates the operation at a call upon a queue is detected, and it uses immediate traffic condition information to determine the control strategy for the subsequent time interval. While the reactive approach allows to remedy the bottleneck with real-time feedback operations, it has limitations related to time lag between the occurrence of congestion and applied control  \cite{khondaker_variable_2015}.
 In contrast, the proactive approach has the capability of acting proactively, while anticipating the behavior of traffic flow \cite{khondaker_variable_2015}.
Thus, it can predict bottleneck formations before they even occur, while potential shock waves can be resolved by restricting traffic inflow. In addition, the nature of predictions of proactive VSL methods allows for a systematic approach for network-wide coordination which supports system optimization, whereas reactive approach is restrained to a localized control logic.


\subsubsection{Reactive Speed Harmonization}

The first field implementation of SPD-HARM was
the VSL system in the German motorway A8 corridor in Munich extented
to the boundary of Salzburg, Austria in 1965 \cite{schick_influence_2003}.
During the early 1960s, US implemented SPD-HARM using VMS on a portion of the New Jersey Turnpike
\cite{robinson_examples_2000}. These SPD-HARM
systems required human interventions to determine the messages or
speed limits based on the conditions such as weather, traffic congestion
and construction schedules. Since 1970s, advances in sensor technologies
and traffic control systems allowed the SPD-HARM to automatically
operate based on the traffic flow or weather conditions using various
types of detectors. The earlier VMS and VSL implementations were often
at the purpose of addressing safety issues under work zone areas or
inclement weather conditions \cite{robinson_examples_2000}. In 2007,
SPD-HARM started focusing on improving
traffic flow mobility. The VSL systems implemented in the M42 motorway at Birmingham,
UK, and Washington State Department of Transportation \cite{brinckerhoff_active_2008},
use algorithms which are automatically activated based on pre-defined
threshold of flow and speed measured by detectors embedded in the
pavement. The systems display the lowered speed limit within a ``control
zone" of a pre-defined length \cite{mott_macdonald_ltd._atm_2008}.

There has been a significant amount of VSL algorithms proposed in the literature to date. 
A reactive approach for VSL to improve
safety and mobility at work-zone areas \cite{park_development_2003,lin_exploring_2004} outperformed the existing VSL algorithms, especially
with the traffic demand fluctuations 
\cite{park_development_2003}. A simulation-based study show that the performance
of VSL is a function of the traffic volume levels \cite{juan_simulation_2004}. 
After reaching a particular
traffic volume level, the benefit can become more apparent, and therefore, VSL needs to
be integrated with ramp metering control \cite{juan_simulation_2004}.
Another VSL algorithm, which implemented
in the Twin Cities Metropolitan area, can identify the moving
jam based on the deceleration rate between adjacent spots \cite{kwon_minnesota_2011}. The field evaluation showed
the reduction in average maximum deceleration by 20\% over the
state-of-the-art in that area while improved the vehicle throughput at the bottleneck
areas. Through its evolution, reactive SPD-HARM has consistently
showed improvements in many aspects such as reliability, safety and
environmental sustainability by providing adequate feedback to the
dynamic traffic conditions. However, the capability of reactive control
is limited as it can be only effective after a bottleneck occurs while it mainly
depends on heuristics.


\subsubsection{Proactive Speed Harmonization}

The necessity of a systematic approach for preventing adverse
impacts from impending shock waves eventually led to the development of the proactive
VSL. The proactive
VSL approach was first proposed in \cite{alessandri_optimal_1998} adopting Kalman Filter aimed at estimating impending traffic status
based on the time-series traffic measurements \cite{welch_introduction_2001}.
Given the estimated traffic flow, the proactive VSL approach derives
a control policy that minimizes various cost functions (e.g.,
average travel time, summation of square densities of all sections). Although this effort initiated prediction-based VSL systems, the prediction using a time-series approach is not robust,
especially under unexpected traffic flow disturbances, since it 
heavily relies on the empirical patterns.

A pioneering effort in developing a proactive VSL system was made
in \cite{hegyi_optimal_2003} using
 model predictive control (MPC). The key aspect
of that work is that it prevents traffic breakdown by decreasing
the density of approaching traffic rather than focusing on reducing
the speed variances. Using MPC, which enabled a network-wide
optimization, a series of VSL systems can be coordinated for system-wide optimization that eventually aims at preventing upstream delays. Another MPC-based proactive VSL system was proposed in \cite{lu_new_2010} focusing on creating a discharge section immediate upstream of
the bottleneck to regulate traffic flow into the bottleneck that remains
close to its capacity. With the intention to influence the motorway mainstream, a traffic flow control approach was proposed in \cite{carlson_optimal_2010}, using a suitable feasible-direction algorithm \cite{papageorgiou_feasible_1995}, that can yield  feedback control policies. These approaches have showed substantial improvements in vehicle
throughput, safety, equity, and driver acceptance through microscopic
simulation studies \cite{carlson_optimal_2010,hegyi_optimal_2003,hegyi_model_2005,hegyi_optimal_2005}.
However, there are significant challenges in practical applications associated with computational requirements. 
Using shock wave theory, another VSL algorithm was proposed in \cite{hegyi_specialist:_2008} that predicts future traffic evolution based on the different traffic states along the freeway segments. By identifying
the location of the front boundaries of shock waves and the active
speed limits, the algorithm maximizes
the discharge rate at the bottleneck \cite{hegyi_specialist:_2008}.

The performance of SPD-HARM can vary depending on the control
approach, characteristics of the topology, and driving behavior. The
potential travel time improvements through SPD-HARM have been debatable during peak
hours \cite{kwon_minnesota_2011,roberts_i-70_2012}. However,
it has been widely agreed that SPD-HARM increases
vehicle throughput at the bottleneck. It has been shown that the vehicle throughput can be increased
by 4-5\% via VSL system \cite{mott_macdonald_ltd._atm_2008}
and by 5-10 \% via rolling SPD-HARM implemented in European
countries with significant benefits in safety since personal injury crashes reduced about 30-35\% \cite{fuhs_synthesis_2010}. The
environmental impacts of SPD-HARM were also substantial demonstrating reduction in vehicle
emissions by 4-10\% (depending on the pollutants) \cite{fuhs_synthesis_2010},
and fuel consumption by 4\% \cite{mott_macdonald_ltd._atm_2008}.

Although previous research efforts reported in the literature have aimed at enhancing our understanding of SPD-HARM algorithms, deriving in real time an optimal solution for each individual vehicle under the hard safety constraints still remains a challenging control problem. In this paper, we address the problem of controlling the speed of a number of automated vehicles before they enter a speed reduction zone on a freeway. We formulate the control problem and provide an analytical, closed-form, optimal solution that can be implemented in real time. The solution yields the optimal acceleration/deceleration for each vehicle in the upstream, and thus it controls the time that each vehicle enters the speed reduction zone.  Furthermore, we provide the conditions under which the rear-end collision avoidance constraint does not become active at any time within the control and speed reduction zones. By controlling the time of reaching the speed reduction zone, the recovery time from the congested area is reduced which, in turn, leads to the increase in average speed, and eventually,  travel time. In addition,  vehicles avoid getting into a stop-and-go driving mode, thereby conserving momentum and energy. The unique contribution of this paper hinges on the following three elements: (1) the formulation of the problem of controlling the speed of a number of automated vehicles before they enter a speed reduction zone on a freeway, (2) a rigorous analysis that reveals the conditions under which the rear-end collision avoidance constraint does not become active, and (3) the proposed optimal control framework bridges the gap between the existing algorithms and the real-world implications by providing an optimal solution that can be implemented in a vehicle in real time.


\subsection{Organization of the Paper}

The remaining paper proceeds as follows. In Section II
we introduce the model, present the assumptions of our
approach and formulate the optimal control problem. In Section III, we provide the control framework, derive an analytical closed-form solution, and discuss the conditions under which the rear-end collision avoidance constraint does not become active. Finally, we provide simulation results and discusion in
Section IV and concluding remarks in Section V.


\section{Problem Formulation}\label{sec:2}
\subsection{Modeling Framework} \label{sec:2a}

\begin{figure}
\centering \includegraphics[width=3.5 in]{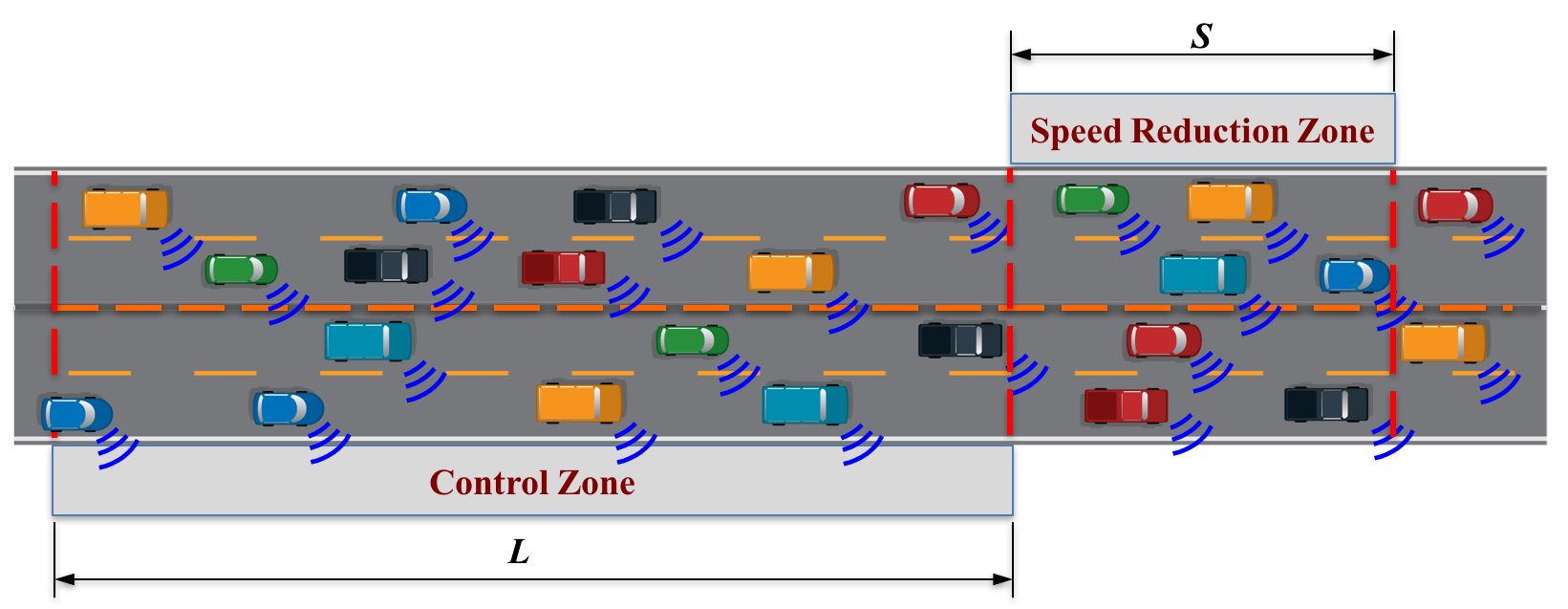} \caption{Automated vehicles within a control zone approaching a speed reduction zone.}
\label{fig:1} 
\end{figure}

We address the problem of optimally controlling the speed of vehicles cruising on a freeway (Fig. \ref{fig:1}) that includes a \textit{speed reduction zone} of a length $S$. The speed reduction zone 
is a bottleneck that builds up as vehicles exceed the bottleneck capacity. 
The freeway has a \textit{control zone} right before the speed reduction zone, inside of which the vehicles need to accelerate/decelerate optimally so as to enter the speed reduction zone with the appropriate speed. Therefore, the speed of the potential queue built-up inside the reduction zone is controlled, and thus the congestion recovery time can be optimized. The distance from the entry of the control zone until the entry of the speed reduction zone is $L$.

We consider a number of automated vehicles
$N(t)\in\mathbb{N}$ in each lane, where $t\in\mathbb{R}^{+}$
is the time, entering the control zone (Fig. \ref{fig:1}). 
Let $\mathcal{N}(t)=\{1,\ldots,N(t)\},$ be the queue in one of the lanes
inside the control zone. The dynamics of each
vehicle $i\in\mathcal{N}(t)$ are represented by
a state equation 
\begin{equation}
\dot{x_{i}}=f(t,x_{i},u_{i}),\qquad x_{i}(t_{i}^{0})=x_{i}^{0},\label{eq:model}
\end{equation}
where $t\in\mathbb{R}^{+}$, $x_{i}(t)$, $u_{i}(t)$
are the state of the vehicle and control input, $t_{i}^{0}$ is the
time that vehicle $i$ enters the control zone, and $x_{i}^{0}$ is
the value of the initial state. For simplicity, we model each
vehicle as a double integrator, e.g., $\dot{p}_{i} =v_{i}(t)$ and $
\dot{v}_{i} =u_{i}(t)$, where $p_{i}(t)\in\mathcal{P}_{i}$, $v_{i}(t)\in\mathcal{V}_{i}$,
and $u_{i}(t)\in\mathcal{U}_{i}$ denote the position, speed and acceleration/deceleration
(control input) of each vehicle $i\in\mathcal{N}(t)$ inside the control zone. Let
$x_{i}(t)=\left[\begin{array}{cc}
p_{i}(t) & v_{i}(t)\end{array}\right]^{T}$ denote the state of each vehicle $i$, with initial value $x_{i}^{0}=\left[\begin{array}{cc}
0 & v_{i}^{0}\end{array}\right]^{T}$, taking values in the state space $\mathcal{X}_{i}=\mathcal{P}_{i}\times\mathcal{V}_{i}$.
The sets $\mathcal{P}_{i}$, $\mathcal{V}_{i}$ and $\mathcal{U}_{i}$,
$i\in\mathcal{N}(t),$ are complete and totally bounded subsets of
$\mathbb{R}$. The state space $\mathcal{X}_{i}$ for each vehicle
$i$ is closed with respect to the induced topology on $\mathcal{P}_{i}\times\mathcal{V}_{i}$
and thus, it is compact.

We need to ensure that for any initial state $(t_{i}^{0},x_{i}^{0})$
and every admissible control $u(t)$, the system \eqref{eq:model}
has a unique solution $x(t)$ on some interval $[t_{i}^{0},t_{i}^{m}]$,
where $t_{i}^{m}$ is the time that vehicle $i\in\mathcal{N}(t)$
enters the speed reduction zone. To ensure that the control input and
vehicle speed are within a given admissible range, the following constraints
are imposed:
\begin{equation}%
\begin{split}
u_{i,min} &  \leqslant u_{i}(t)\leqslant u_{i,max},\quad\text{and}\\
0 &  \leqslant v_{min}\leqslant v_{i}(t)\leqslant v_{max},\quad\forall
t\in\lbrack t_{i}^{0},t_{i}^{m}],
\end{split}
\label{speed_accel constraints}%
\end{equation}
where $u_{i,min}$, $u_{i,max}$ are the minimum and maximum
control inputs (maximum deceleration/ acceleration) for each vehicle $i\in\mathcal{N}(t)$, and $v_{min}$, $v_{max}$
are the minimum and maximum speed limits respectively. For simplicity, in the
rest of the paper we do not consider vehicle diversity, and thus, we set $u_{i,min}=u_{min}$
and $u_{i,max}=u_{max}$.

To ensure the absence of any rear-end collision of two consecutive vehicles
traveling on the same lane, the position of the preceding vehicle should be
greater than or equal to the position of the following vehicle plus a
 safe distance $\delta(t)<S$, which is a function of speed. Thus, we impose the rear-end safety
constraint
\begin{equation}
s_{i}(t)=p_{i-1}(t)-p_{i}(t)\geqslant\delta(v_{ave}(t)),~\forall t\in\lbrack t_{i}^{0}%
,t_{i}^{m}],\label{eq:rearend}%
\end{equation}
where $i-1$ denotes the vehicle which is physically
immediately ahead of $i$ in the same lane, and $v_{ave}(t)$ is the average speed of the vehicles inside the control zone at time $t$.

In the modeling framework described above, we impose the following
assumptions: 


\begin{customass}{1} \label{ass:lane} The vehicles do not change lanes inside the control and speed reduction zones.
\end{customass}


\begin{customass}{2} \label{ass:srz} Each vehicle cruises inside the speed reduction zone with the imposed speed limit, $v_{srz}$.
\end{customass}


\begin{customass}{3} \label{ass:sensor} Each vehicle $i$ has proximity
sensors and can measure local information without errors or delays.
\end{customass}

We briefly comment on the above assumptions. The first assumption is considered to simplify the problem and limit the scope of the paper on understanding the implications of the proposed approach on each lane seperately. 
The second assumption 
is intended to enhance safety awareness, but it could be modified appropriately, if necessary.
The third assumption may be a strong requirement to impose. 
However, it is relatively straightforward to extend our
framework in the case that it is relaxed, as long as  the noise
in the measurements and/or delays are bounded. For example, we can
determine the uncertainties of the state of the vehicle as a result
of sensing or communication errors and delays, and incorporate these
into the safety constraints.


\subsection{Optimal Control Problem Formulation}
\label{sec:2b}

We consider the problem of deriving the optimal acceleration/deceleration
of each automated vehicle in a freeway (Fig. \ref{fig:1}), under the hard safety
constraint to avoid rear-end collision. The potential benefits of
the solution of this problem are substantial. By controlling the speed of the vehicles
in the upstream or tighten the inflow traffic, the speed of queue
built-up decreases, and thus the congestion recovery time is also reduced.
Even though the speed of each vehicle is reduced, the throughput at the speed reduction zone is maximized.
Moreover, by minimizibg the acceleration/deceleration of each vehicle,
we minimize transient engine operation, thus we can have direct benefits in
fuel consumption \cite{Malikopoulos2008b} and emissions since internal combustion
engines are optimized over steady state operating points (constant torque and
speed) \cite{Malikopoulos2010a,Malikopoulos2015a}.

When a vehicle $i$ enters the control zone, it receives some information from the vehicle $i-1$ which is 
physically located ahead of it.

\begin{customdef}{1}\label{def:1} For each vehicle
$i$ when it enters a control zone, we define the \textit{information
set} $Y_{i}(t)$ as 
\begin{equation}
Y_{i}(t)\triangleq\Big\{ p_{i}(t),v_{i}(t),t_{i}^{m}\Big\},~\forall t\in\lbrack t_{i}^{0},t_{i}^{m}],\label{eq:local}
\end{equation}
\end{customdef} where $p_{i}(t),v_{i}(t)$ are the position and speed
of vehicle $i$ inside the control zone, and $t_{i}^{m}$,
is the time targeted for vehicle $i$ to enter the
speed reduction zone. Each vehicle $i\in\mathcal{N}(t)$ has proximity sensors and can observe and/or
estimate the information in $Y_{i}(t)$ without errors or delays (Assumption
\ref{ass:sensor}). Note that once vehicle $i$ enters the control zone, then immediately all information in
$Y_{i}(t)$ becomes available to $i$: $p_{i}(t),v_{i}(t)$ are read from its
sensors and $t_{i}^{m}$ can also be computed
at that time based on the information the vehicle $i$ receives from $i-1$ as described next. 

The time $t_{i}^{m}$ that the vehicle $i$ will be entering the speed reduction zone is
restricted by the imposing rear-end collision constraint. Therefore, to ensure that \eqref{eq:rearend} is satisfied at $t_{i}^{m}$ we impose the following condition
\begin{equation}
t_{i}^{m}=\max\Bigg\{ \min\Big\{t_{i-1}^{m}+\frac{\delta(v_{ave}(t))}{v_{i-1}(t_{i-1}^{m})},\frac{L}{v_{min}}\Big\}, \frac{L}{v_{i}(t_{i}^{0})}, \frac{L}{v_{max}}\Bigg\},\label{eq:lem1b}%
\end{equation}
where $v_{i-1}(t_{i-1}^{m})$ is the speed of the vehicle $i-1$ at the time $t_{i-1}^{m}$  that enters the speed reduction zone, and it is equal to the speed, $v_{srz}$, imposed inside the reduction zone (Assumption \ref{ass:srz}). Thus, the condition \eqref{eq:lem1b} ensures that the time $t_{i}^{m}$ that vehicle $i$ will be entering the speed reduction zone is feasible and can be attained based on the imposed speed limits inside the control zone. In addition, for low traffic flow where vehicle $i-1$ and $i$ might be located far away from each other, there is no compelling reason for vehicle $i$ to accelerate within the control zone just to have a distance $\delta(v_{ave}(t))$ from vehicle $i-1$ at the time $t_{i}^{m}$ that vehicle $i$ enters the speed reduction zone. Therefore, in such cases vehicle $i$ can keep cruising within the control zone with the initial speed $v_{i}(t_{i}^{0})$ that entered the control zone at $t_{i}^{0}.$ The recursion is initialized when the first vehicle enters
the control zone, i.e., it is assigned $i=1$. In this case, $t_{1}^{m}$ can be
externally assigned as the desired exit time of this vehicle whose behavior is
unconstrained. Thus the time $t_{1}^{m}$ is 
fixed and available through $Y_{1}(t)$. The
second vehicle will access $Y_{1}(t)$ to compute the times $t_{2}^{m}$. The third vehicle will access $Y_{2}(t)$ and the
communication process will continue with the same fashion until the
vehicle $N(t)$ in the queue access the $Y_{N(t) -1}(t)$.

We consider the problem of minimizing the control input (acceleration/deceleration) for each vehicle from the time $t_{i}^{0}$ that the vehicle $i$ enters the control
zone until the time $t_{i}^{m}$ that it enters the speed reduction zone
under the hard safety constraint to avoid rear-end collision. 
By minimizing each vehicle's acceleration/deceleration, we minimize transient engine
operation \cite{Malikopoulos2015b, malikopoulos2013stochastic}, and thus, we can have direct benefits in fuel consumption  and
emissions since internal combustion engines are optimized over steady state
operating points (constant torque and speed).

The
optimal control problem is formulated as to minimize the L$^2$ norm of the control input 
\begin{gather}
\min_{u_{i}}\frac{1}{2}\int_{t_{i}^{0}}^{t_{i}^{m}}u_{i}^{2}%
(t)~dt,\label{eq:decentral}\\
\text{subject to}:\eqref{eq:model}~\text{and}~ \eqref{speed_accel constraints},\nonumber
\end{gather}
with initial and final conditions: $p_{i}(t_{i}^{0})=0$, $p_{i}(t_{i}^{m})=L,$ $t_{i}^{0}$, $v_{i}(t_{i}^{0}),$ $t_{i}^{m},$ and $v_{i}(t_{i}^{m})=v_{srz}$. Note that we have omitted the rear end safety constraint \eqref{eq:rearend} in the problem formulation above. 
The analytical solution of the problem including the rear-end collision avoidance
constraint may become intractable for real-time implementation.  
However, in the following section, we provide the conditions under which the rear-end collision avoidance constraint does not become active
at any time in $(t_{i}^{0},t_{i}^{m})$ as long as it is not active
at $t=t_{i}^{0}$. Note that we can guarantee rear-end collision
avoidance at time $t_{i}^{m}$ based on \eqref{eq:lem1b}.


\section{Solution of the optimal control problem}

\label{sec:3}

For the analytical solution and real-time implementation of the control
problem \eqref{eq:decentral}, we apply Hamiltonian analysis. In our
analysis, we consider that when the
vehicles enter the control zone, none of the constraints are active.
However, this might not be in general true. For example, a vehicle
may enter the control zone with speed higher than the speed limit.
In this case, we need to solve an optimal control problem starting
from an infeasible state. 


The Hamiltonian function is formulated as follows 
\begin{gather}
H_{i}\big(t,x(t),u(t)\big)=L_{i}\big(t,x(t),u(t)\big)+\lambda^{T}\cdot f_{i}\big(t,x(t),u(t)\big)\\
+\mu^{T}\cdot g_{i}\big(t,x(t),u(t)\big),\nonumber 
\end{gather}
where 
\[
g_{i}\big(t,x(t),u(t)\big)=\left\{ \begin{array}{ll}
u_{i}(t)-u_{max}\leq0,\\
u_{min}-u_{i}(t)\leq0,\\
v_{i}(t)-v_{max}\leq0,\\
v_{min}-v_{i}(t)\leq0.
\end{array}\right.
\]
From \eqref{eq:decentral}, the state equations \eqref{eq:model}, and the
control/state constraints \eqref{speed_accel constraints}, for each vehicle
$i\in\mathcal{N}(t)$ the Hamiltonian function becomes
\begin{gather}
H_{i}\big(t,p(t),v(t),u(t)\big)=\frac{1}{2}u_{i}^{2}+\lambda_{i}^{p}\cdot
v_{i}+\lambda_{i}^{v}\cdot u_{i}\nonumber\\
+\mu_{i}^{a}\cdot(u_{i}-u_{max})+\mu_{i}^{b}\cdot(u_{min}-u_{i})+\mu_{i}%
^{c}\cdot(v_{i}-v_{max})\nonumber\\
+\mu_{i}^{d}\cdot(v_{min}-v_{i}),\label{eq:16b}%
\end{gather}
where $\lambda_{i}^{p}$ and $\lambda_{i}^{v}$ are the costates, and $\mu^{T}$
is a vector of Lagrange multipliers with
\begin{equation}
\mu_{i}^{a}=\left\{
\begin{array}
[c]{ll}%
>0, & \mbox{$u_{i}(t) - u_{max} =0$},\\
=0, & \mbox{$u_{i}(t) - u_{max} <0$},
\end{array}
\right.  \label{eq:17a}%
\end{equation}%
\begin{equation}
\mu_{i}^{b}=\left\{
\begin{array}
[c]{ll}%
>0, & \mbox{$u_{min} - u_{i}(t) =0$},\\
=0, & \mbox{$u_{min} - u_{i}(t)<0$},
\end{array}
\right.  \label{eq:17b}%
\end{equation}%
\begin{equation}
\mu_{i}^{c}=\left\{
\begin{array}
[c]{ll}%
>0, & \mbox{$v_{i}(t) - v_{max} =0$},\\
=0, & \mbox{$v_{i}(t) - v_{max}<0$},
\end{array}
\right.  \label{eq:17c}%
\end{equation}%
\begin{equation}
\mu_{i}^{d}=\left\{
\begin{array}
[c]{ll}%
>0, & \mbox{$v_{min} - v_{i}(t)=0$},\\
=0, & \mbox{$v_{min} - v_{i}(t)<0$}.
\end{array}
\right.  \label{eq:17d}%
\end{equation}
The Euler-Lagrange equations become
\begin{equation}
\dot{\lambda}_{i}^{p}=-\frac{\partial H_{i}}{\partial p_{i}}=0,\label{eq:EL1}%
\end{equation}
and
\begin{equation}
\dot{\lambda}_{i}^{v}=-\frac{\partial H_{i}}{\partial v_{i}}=\left\{
\begin{array}
[c]{ll}%
-\lambda_{i}^{p}, & \mbox{$v_{i}(t) - v_{max} <0$}~\text{and}\\
& \mbox{$v_{min} - v_{i}(t)<0$},\\
-\lambda_{i}^{p}-\mu_{i}^{c}, & \mbox{$v_{i}(t) - v_{max} =0$},\\
-\lambda_{i}^{p}+\mu_{i}^{d}, & \mbox{$v_{min} - v_{i}(t)=0$},
\end{array}
\right.  \label{eq:EL2}%
\end{equation}
with given initial and final conditions $p_{i}(t_{i}^{0})=0$, $p_{i}(t_{i}^{m})=L$, $v_{i}(t_{i}^{0})$, and $v_{i}(t_{i}^{m})$.
The necessary condition for optimality is
\begin{equation}
\frac{\partial H_{i}}{\partial u_{i}}=u_{i}+\lambda_{i}^{v}+\mu_{i}^{a}%
-\mu_{i}^{b}=0. \label{eq:KKT1}%
\end{equation}

To address this problem, the constrained and unconstrained arcs need to
be pieced together to satisfy the Euler-Lagrange equations and necessary
condition of optimality. The analytical solution of \eqref{eq:decentral}
without considering state and control constraints was presented in earlier papers
 \cite{Rios-Torres2015,Rios-Torres2,Ntousakis:2016aa} for coordinating in real time CAVs at highway on-ramps and \cite{ZhangMalikopoulosCassandras2016} at two adjacent intersections. 

When the state and control constraints are not active, we
have $\mu_{i}^{a}=\mu_{i}^{b}=\mu_{i}^{c}=\mu_{i}^{d}=0.$ Applying
the necessary condition \eqref{eq:KKT1}, the optimal control can
be given $u_{i}+\lambda_{i}^{v}=0$,$\quad i\in\mathcal{N}(t).$ The Euler-Lagrange equations yield 
$\dot{\lambda}_{i}^{p}=-\frac{\partial H_{i}}{\partial p_{i}}=0$ and $\dot{\lambda}_{i}^{v}=-\frac{\partial H_{i}}{\partial v_{i}}=-\lambda_{i}^{p}.$ From the former equation we have $\lambda_{i}^{p}=a_{i}$ and from the latter  
$\lambda_{i}^{v}=-(a_{i}t+b_{i})$, where $a_{i}$ and $b_{i}$
are constants of integration corresponding to each vehicle $i$. Consequently,
the optimal control input (acceleration/deceleration) as a function
of time is given by 
\begin{equation}
u_{i}^{*}(t)=a_{i}t+b_{i}.\label{eq:20}
\end{equation}
Substituting the last equation into the vehicle dynamics equations
\eqref{eq:model}, we can find the optimal speed and position for
each vehicle, namely 
\begin{equation}
v_{i}^{*}(t)=\frac{1}{2}a_{i}t^{2}+b_{i}t+c_{i}\label{eq:21}
\end{equation}
\begin{equation}
p_{i}^{*}(t)=\frac{1}{6}a_{i}t^{3}+\frac{1}{2}b_{i}t^{2}+c_{i}t+d_{i},\label{eq:22}
\end{equation}
where $c_{i}$ and $d_{i}$ are constants of integration. These constants
can be computed by using the initial and final conditions. Since we
seek to derive the optimal control (\ref{eq:20}) in real time, we can designate
initial values $p_{i}(t_{i}^{0})$ and $v_{i}(t_{i}^{0})$, and initial
time, $t_{i}^{0}$, to be the current values of the states $p_{i}(t)$
and $v_{i}(t)$ and time $t$, where $t_{i}^{0}\le t\le t_{i}^{m}$.
Therefore the constants of integration will be functions of time and
states, i.e., $a_{i}(t,p_{i},v_{i}),b_{i}(t,p_{i},v_{i}),c_{i}(t,p_{i},v_{i})$,
and $d_{i}(t,p_{i},v_{i})$. It follows that
\eqref{eq:21} and \eqref{eq:22}, along with the initial and terminal conditions,
can be used to form a system of four equations of the form $\mathbf{T}%
_{i}\mathbf{b}_{i}=\mathbf{q}_{i}$, namely%

\begin{equation}
\left[
\begin{array}
[c]{cccc}%
\frac{1}{6}t^{3} & \frac{1}{2}t^{2} & t & 1\\
\frac{1}{2}t^{2} & t & 1 & 0\\
\frac{1}{6}(t_{i}^{m})^{3} & \frac{1}{2}(t_{i}^{m})^{2} & t_{i}^{m} & 1\\
\frac{1}{2}(t_{i}^{m})^{2} & t_{i}^{m} & 1 & 0
\end{array}
\right]  .\left[
\begin{array}
[c]{c}%
a_{i}\\
b_{i}\\
c_{i}\\
d_{i}%
\end{array}
\right]  =\left[
\begin{array}
[c]{c}%
p_{i}(t)\\
v_{i}(t)\\
p_{i}(t_{i}^{m})\\
v_{i}(t_{i}^{m})
\end{array}
\right], \label{eq:23}%
\end{equation}
where $t_{i}^{m}$ is specified by \eqref{eq:lem1b}, and $v_{i}(t_{i}^{m})=v_{srz}$ is the imposed speed limit at the speed reduction zone.
Hence, we have
\begin{equation}
\mathbf{b}_{i}(t,p_{i}(t),v_{i}(t))=(\mathbf{T}_{i})^{-1}.\mathbf{q}%
_{i}(t,p_{i}(t),v_{i}(t)), \label{eq:24}%
\end{equation}
where $\mathbf{b}_{i}(t,p_{i}(t),v_{i}(t))$ contains the four integration
constants $a_{i}(t,p_{i},v_{i})$, $b_{i}(t,p_{i},v_{i})$, $c_{i}(t,p_{i}%
,v_{i})$, $d_{i}(t,p_{i},v_{i})$. Thus, \eqref{eq:20} can be written as
\begin{equation}
u_{i}^{\ast}(t,p_{i}(t),v_{i}(t))=a_{i}(t,p_{i}(t),v_{i}(t))t+b_{i}%
(t,p_{i}(t),v_{i}(t)). \label{eq:25}%
\end{equation}
Since \eqref{eq:23} can be computed in real time, the controller yields the
optimal control in real time for each vehicle $i$, with feedback provided
through the re-calculation of the vector $\mathbf{b}_{i}(t,p_{i}(t),v_{i}(t))$
in \eqref{eq:24}. Similar results can be obtained when
the state and control constraints become active \cite{Malikopoulos2017}.

This analytical solution for the case when the state and control constraints are not active, however, does not include the rear-end collision avoidance constraint.
Thus, we investigate the conditions under which the rear-end collision avoidance constraint does not become active for any two vehicles $i-1$ and $i$ at any time in $[t_{i}^{0},t_{i}^{m}]$, if they are not active at $t=t_{i}^{0}$.

\begin{customthm}{1}\label{theo:1}
Suppose that
there exists a feasible solution of the control problem \eqref{eq:decentral}. Then, if
$p_{i-1}(t_{i}^{0})-p_{i}(t_{i}^{0})=l\ge\delta(v_{ave}(t))$ for $i=2,3,\ldots$, the
constraint \eqref{eq:rearend} does not become active, i.e., $p_{i-1}(t)-p_{i}(t)\geq\delta(v_{ave}(t))$ for
all $t\in(t_{i}^{0},t_{i}^{m}]$, if 
\begin{gather}%
\big(v_{i-1}(t_{i}^{0})- v_{i}(t_{i}^{0})\big)\\
 \ge \frac{\big(l - \delta(v_{ave}(t))\big)\cdot \big( 2~t^3-3~t^2~t^{m}_{i}-(t^{m}_{i})^3 \big)}{t^{m}_{i} \cdot (t-t^{m}_{i})^2 \cdot t}.
\label{theo:1aa}
\end{gather}
\end{customthm}

\begin{proof}
Without loss of generality and to simplify notation, we reset the time at $t=t^0_{i}$, i.e., $t^0_{i}=0$. Thus, at $t^0_{i}=0$, vehicle $i-1$ has traveled a distance $l$ inside the control zone and has a speed $v_{i-1}(t^0_{i})$. Similarly, at $t^0_{i}=0$, vehicle $i$ just entered the control zone with an initial speed $v_{i}(t_{i}^{0})=v_{i}^0$. From \eqref{eq:lem1b}, the position of vehicle $i-1$ will be $L+\delta(v_{ave}(t))$ at $t=t^{m}_{i}$ and will have the speed imposed at the speed reduction zone, $v_{i-1}(t^{m}_{i})=v_{srz}$. Similarly, at $t=t^{m}_{i}$ the vehicle $i$ will be at the entry of the speed reduction zone, so its position will be $L$ and its speed will be the speed imposed at the speed reduction zone, $v_{i}(t^{m}_{i})=v_{srz}$. 

Since the state and control constraints are not active, the control input, speed and position for each vehicle $i$ are given by \eqref{eq:20} - \eqref{eq:22}.
Substituting the conditions at $t^0_{i}=0$ and $t=t^{m}_{i}$ for the vehicles $i-1$ and $i$ in \eqref{eq:23}, and by letting $\sigma_{i}=\frac{1}{6}(t^{m}_{i})^{3}$, $\rho_{i}=\frac{1}{2}(t^{m}_{i})^{2}$, and $\tau_{i}=t^{m}_{i}$ we have
\begin{equation}
\left[
\begin{array}
[c]{cccc}%
0 & 0 & 0 & 1\\
0 & 0 & 1 & 0\\
\sigma_{i} & \rho_{i} & \tau_{i} & 1\\
\rho_{i} & \tau_{i} & 1 & 0
\end{array}
\right] . \left[
\begin{array}
[c]{c}%
a_{i-1}\\
b_{i-1}\\
c_{i-1}\\
d_{i-1}%
\end{array}
\right] =\left[
\begin{array}
[c]{c}%
l\\
v_{i-1}(0)\\
L+\delta(v_{ave}(t))\\
v_{i-1}(t^{m}_{i})
\end{array}
\right] ,\label{theo:1a}%
\end{equation}
and
\begin{equation}
\left[
\begin{array}
[c]{cccc}%
0 & 0 & 0 & 1\\
0 & 0 & 1 & 0\\
\sigma_{i} & \rho_{i} & \tau_{i} & 1\\
\rho_{i} & \tau_{i} & 1 & 0
\end{array}
\right] . \left[
\begin{array}
[c]{c}%
a_{i}\\
b_{i}\\
c_{i}\\
d_{i}%
\end{array}
\right] =\left[
\begin{array}
[c]{c}%
0\\
v_{i}^0\\
L\\
v_{i}(t^{m}_{i})
\end{array}
\right], \label{theo:1b}%
\end{equation}
where $v_{i-1}(0)$ is the speed of vehicle $i-1$ at time $t^0_{i}=0$, $v_{i-1}(t^{m}_{i}) = v_{i}(t^{m}_{i})$$=v_{srz}$ is the speed at the reduction zone; and $a_{k}, b_{k}, c_{k}$, and $d_{k}$, $k=i-1,i,$ are the constants of integration that can be computed from \eqref{eq:24} where
\begin{equation}%
(\mathbf{T}_{i})^{-1} = \frac{1}{\text{det}({T}_{i})} \cdot \text{adj}(T), 
\label{theo:1c}
\end{equation}
and $\text{det}({T}_{i})=-\sigma_{i}\cdot\tau_{i}+\rho_{i}^2=\frac{(t^{m}_{i})^4}{12}.$

Thus from \eqref{theo:1a} and \eqref{theo:1b} and by letting $\gamma_{i}= \sigma_{i}-\rho_{i}\cdot \tau_{i}$ we have
\begin{gather}
\left[
\begin{array}
[c]{c}%
a_{i-1}\\
b_{i-1}\\
c_{i-1}\\
d_{i-1}%
\end{array}
\right] = \frac{1}{\text{det}({T}_{i})} \\
\cdot \left[
\begin{array}
[c]{cccc}%
\tau_{i} & -\rho_{i}+\tau_{i}^2 & -\tau_{i} & \rho_{i}\\
-\rho_{i} & \gamma_{i} & \rho_{i} & -\sigma_{i}\\
0 & \text{det}({T}_{i}) & 0 & 0\\
\text{det}({T}_{i}) & 0 & 0 & 0
\end{array}
\right]\cdot  
\left[
\begin{array}
[c]{c}%
l\\
v_{i-1}(0)\\
L+\delta(v_{ave}(t))\\
v_{i-1}(t^{m}_{i})
\end{array}
\right] 
\label{theo:1c}
\end{gather}
and
\begin{gather}\label{theo:1d}
\left[
\begin{array}
[c]{c}%
a_{i}\\
b_{i}\\
c_{i}\\
d_{i}%
\end{array}
\right] = \frac{1}{\text{det}({T}_{i-1})} \\
\cdot  
\left[
\begin{array}
[c]{cccc}%
\tau_{i} & -\rho_{i}+\tau_{i}^2 & -\tau_{i} & \rho_{i}\\
-\rho_{i} & \gamma_{i} & \rho_{i} & -\sigma_{i}\\
0 & \text{det}({T}_{i}) & 0 & 0\\
\text{det}({T}_{i}) & 0 & 0 & 0
\end{array} 
\right] \cdot  
\left[
\begin{array}
[c]{c}%
0\\
v_{i}^0\\
L\\
v_{i-1}(t^{m}_{i})
\end{array}
\right]. 
\end{gather}

From \eqref{theo:1c} and \eqref{theo:1d}, we can compute the constants of integration $a_i,b_i$ 
\begin{gather}%
a_{i-1} = \frac{6}{(t^{m}_{i})^2}\cdot \big(v_{i-1}(0)+ v_{i-1}(t^{m}_{i}) \big)  \\
-\frac{12}{(t^{m}_{i})^3} \cdot(L+\delta(v_{ave}(t))-l), \nonumber \\
a_{i} = \frac{6}{(t^{m}_{i})^2}\cdot \big(v_{i}^{0}+v_{i-1}(t^{m}_{i})\big) -\frac{12}{(t^{m}_{i})^3} \cdot L,
\label{theo:1e}
\end{gather}
\begin{gather}%
b_{i-1} = \frac{6}{(t^{m}_{i})^2}\cdot (L+\delta(v_{ave}(t))-l) \nonumber\\
-\frac{12}{t^{m}_{i}} \cdot \big(\frac{1}{3} v_{i-1}(0) + \frac{1}{6} v_{i-1}(t^{m}_{i}) \big), \\
 b_{i}  = \frac{6}{(t^{m}_{i})^2}\cdot L -\frac{12}{t^{m}_{i}} \cdot \big( \frac{1}{3} v_{i}^0 + \frac{1}{6} v_{i-1}(t^{m}_{i}) \big).
\label{theo:1f}
\end{gather}
The constants of integration $c_i$ and $d_i$, can be derived from \eqref{theo:1a} and \eqref{theo:1b} directly 
\begin{equation}%
c_{i-1}=  v_{i-1}(0), ~c_{i}= v_{i}^{0}, ~d_{i-1}=l, ~\text{and} ~d_{i}=0.
\label{theo:1ff}
\end{equation}

We know $p_{i-1}(t_{i}^{0})-p_{i}(t_{i}^{0})=l\ge\delta(v_{ave}(t))$ for $i=2,3,\ldots$ For all $t\in(t_{i}^{0},t_{i}^{m}]$, we have
\begin{gather}%
p^{*}_{i-1}(t) - p^{*}_{i}(t) \ge \delta(v_{ave}(t))  \Rightarrow \\
\frac{1}{6} (a_{i-1}-a_{i})\cdot t^{3} + \frac{1}{2} (b_{i-1}-b_{i}) \cdot t^{2} + (c_{i-1}-c_{i}) \cdot t \nonumber\\
+ (d_{i-1}-d_{i}) \ge\delta(v_{ave}(t)),~ t\in(t_{i}^{0},t_{i}^{m}].
\label{theo:1g}
\end{gather}
Substituting the constants of integration, \eqref{theo:1g} becomes
\begin{gather}%
\frac{t^3}{6} \Big[ \frac{6}{(t^{m}_{i})^2}\cdot \big(v_{i-1}(0)- v_{i}^0\big) +  \frac{12}{(t^{m}_{i})^3}\cdot \big(\delta(v_{ave}(t))-l\big) \Big]  \nonumber\\
+\frac{t^2}{2} \Big[- \frac{4}{t^{m}_{i}}\cdot \big(v_{i-1}(0)- v_{i}^0\big) -  \frac{6}{(t^{m}_{i})^2} \cdot \big(\delta(v_{ave}(t))-l\big)\Big]\nonumber \\
+ t\cdot \big(v_{i-1}(0)- v_{i}^0\big) - (\delta(v_{ave}(t))-l)\ge 0,~
t\in(t_{i}^{0},t_{i}^{m}]. \nonumber
\label{theo:1h}
\end{gather}
By rearranging the terms in the last equation we have
\begin{gather}%
\big(v_{i-1}(0)- v_{i}^0\big) \cdot\Big(\frac{t^3}{(t^{m}_{i})^2}-\frac{2~t^2}{t^{m}_{i}} +t \Big) \nonumber \\
+\big(\delta(v_{ave}(t))-l\big) \cdot \Big(\frac{2~t^3}{(t^{m}_{i})^3}-\frac{3~t^2}{(t^{m}_{i})^2}-1 \Big) \ge 0, ~t\in(t_{i}^{0},t_{i}^{m}].
\label{theo:1i}
\end{gather}
Since $t^{m}_{i}>0$ we can multiply both sides by $(t^{m}_{i})^3$, hence
\begin{gather}%
\big(v_{i-1}(0)- v_{i}^0\big) \cdot t \cdot t^{m}_{i} \cdot (t-t^{m}_{i})^2 \nonumber\\
+\big(\delta(v_{ave}(t))-l\big)\cdot \big( 2~t^3-3~t^2~t^{m}_{i}-(t^{m}_{i})^3 \big)\ge0, 
\label{theo:1j}
\end{gather}
or
\begin{gather}%
\big(v_{i-1}(t_{i}^{0})- v_{i}(t_{i}^{0})\big) \ge \nonumber\\
 \frac{\big(l- \delta(v_{ave}(t))\big)\cdot \big( 2~t^3-3~t^2~t^{m}_{i}-(t^{m}_{i})^3 \big)}{t^{m}_{i} \cdot (t-t^{m}_{i})^2 \cdot t}, ~t\in(t_{i}^{0},t_{i}^{m}].
\label{theo:1i}
\end{gather}
\end{proof}

\begin{customrem}{1}\label{remark:1}
The function  $f(t)=\frac{\big(l-\delta(v_{ave}(t))\big)\cdot \big( 2~t^3-3~t^2~t^{m}_{i}-(t^{m}_{i})^3 \big)}{t^{m}_{i} \cdot (t-t^{m}_{i})^2 ~ t}, ~t\in(t_{i}^{0},t_{i}^{m}],$ in \eqref{theo:1i} is increasing with respect to $t$ for all $t\in(t_{i}^{0},t_{i}^{m}]$, since $f'(t)=\frac{(t^{m}_{i})^4 +t^3~ t^{m}_{i} -3~ t^2~ (t^{m}_{i})^2 -3~ t~(t^{m}_{i})^3}{(t-t^{m}_{i})^3 ~ t^2}\ge 0$ for all $t\in(t_{i}^{0}, t_{i}^{s}]$, where $t_{i}^{s}$ is some time such that $t_{i}^{0}\le t_{i}^{s} \le t_{i}^{m}$, and it is decreasing with respect to $t$ for all $t\in (t_{i}^{s}, t_{i}^{m}]$. Therefore, it is sufficient to check the condition \eqref{theo:1i} at  $t=t_{i}^{s}$ where $f(t)$ takes its maximum value.
\end{customrem}

\begin{customrem}{2}\label{remark:2}
If at time $t_{i}^{0}$ vehicle $i$ enters the control zone with an initial speed $v_{i}(t_{i}^{0})$ such that there is $t=t_{i}^{s}\in(t_{i}^{0},t_{i}^{m}]$ at which the condition \eqref{theo:1i} does not hold, then vehicle $i$ must update $t_{i}^{m}$ given by \eqref{eq:lem1b}. Let $(v_{i-1}(t_{i}^{0})- v_{i}(t_{i}^{0}) =\zeta$ and suppose that at $t=t_{i}^{s}\in(t_{i}^{0},t_{i}^{m}]$, $\frac{\big(l- \delta(v_{i}(t_{i}^{s}))\big)\cdot \big( 2~(t_{i}^{s})^3-3~(t_{i}^{s})^2~t^{m}_{i}-(t^{m}_{i})^3 \big)}{t^{m}_{i} \cdot (t_{i}^{s}-t^{m}_{i})^2 \cdot t_{i}^{s}}>\zeta$. Then from \eqref{theo:1i} we have
\begin{gather}
(t^{m}_{i})^3 \cdot(t_{i}^{s}~\zeta +l-\delta(v_{ave}(t))) - (t^{m}_{i})^2 ~2~(t^{s}_{i})^2~\zeta \nonumber\\
+ t^{m}_{i}\cdot(3~(t^{s}_{i})^2~(l-\delta(v_{ave}(t)))+(t^{s}_{i})^3~\zeta)\nonumber \\ - 2~(l-\delta(v_{ave}(t)))~(t^{s}_{i})^3 )\ge 0.
\label{rem:2}
\end{gather}
Therefore, the condition \eqref{theo:1i} can be satisfied, if vehicle $i$ enters the speed reduction zone at time $\overline{t^{m}_{i}}$ which is the solution of \eqref{rem:2}.
\end{customrem}

Next, we investigate whether \eqref{eq:rearend} becomes active for any $t\in(t_{i}^{m},t_{i-1}^{f}]$, namely, from the time $t_{i}^{m}$ that vehicle $i$ enters the speed reduction zone until the time $t_{i-1}^{f}$ vehicle $i-1$ exits the speed reduction zone. 

\begin{customthm}{2}\label{theo:2}
Suppose Assumption \ref{ass:srz} is in effect. 
Then, \eqref{eq:rearend} does not become active, when the state and control constraints are not active, for all $t\in(t_{i}^{m},t_{i-1}^{f}]$, if $p_{i-1}(t_{i}^{m})-p_{i}(t_{i}^{m})=l\ge\delta(v_{ave}(t))$.
\end{customthm}

\begin{proof}
Since $v^{*}_{i}(t) = v^{*}_{i-1}(t) = v_{srz}$, for all $t\in(t_{i}^{m},t_{i-1}^{f}]$, then $p_{i-1}(t)-p_{i}(t)=l\ge\delta(v_{ave}(t))$ for all $t\in(t_{i}^{m},t_{i-1}^{f}]$.
\end{proof}

\begin{customcor}{1} \label{corol:1}
The condition  \eqref{theo:1i}  guarantees that \eqref{eq:rearend} does not become active for all $t\in(t_{i}^{0},t_{i-1}^{f}]$ for the case when the state and control constraints are not active.
\end{customcor}

\begin{customrem}{3}\label{remark:3}
If $p_{i-1}(t_{i}^{0})-p_{i}(t_{i}^{0})=l<\delta(v_{ave}(t))$ for $i=2,3,\ldots$, then the
constraint \eqref{eq:rearend} is active at $t = t_{i}^{0}$. In this case, the optimal solution would need to start from a non-feasible state. In such instance, the vehicle must be enforced to decelerate as needed until it reaches the minimum safe distance $\delta(v_{ave}(t))$ from its preceding vehicle before the optimal control applies.
\end{customrem}


\section{Simulation Framework and Results}

\label{sec:4}
\begin{figure}
\centering
\includegraphics[width=3in]{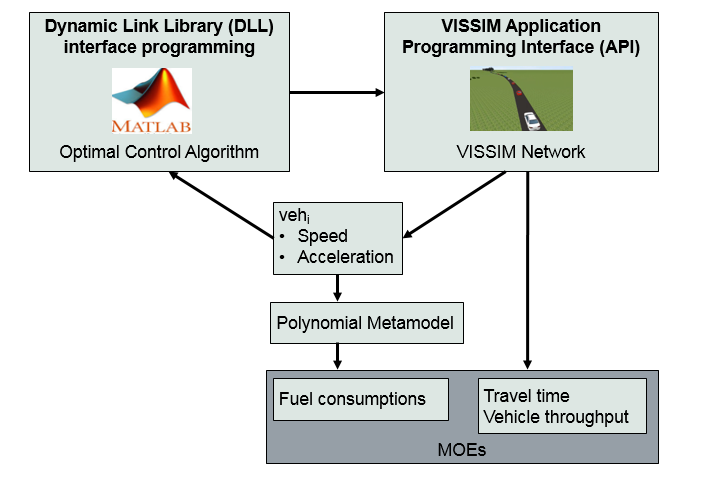}

\caption{Overview of simulation framework}

\label{fig:2}
\end{figure}

To evaluate the effectiveness of our approach, a simulation framework was established  as shown in Fig. \ref{fig:2}.
The analytical, closed-form solution of the proposed control algorithm described
in the previous section was implemented for each vehicle using MATLAB Dynamic
Link Library (DLL) interface programming to allow data exchange with
 external programs within the framework. A simulation test-bed
network was developed under a microscopic traffic simulation software, VISSIM, and it was integrated into the
framework through COM interface. 

The mobility measures such as travel time, average speed and vehicle
throughput were directly obtained from VISSIM. Fuel consumption is a function of vehicle speed and acceleration. The speed of the vehicle is deemed characteristic of the driver's preference, and thus, it might not be effectively controlled as the driver would override any other than the preferred speed. Moreover, if someone considers to minimize fuel consumption with respect to speed, then the question is how to determine a low bound of the speed since the minimum fuel consumption corresponds to zero vehicle speed. Therefore, through the proposed approach, we consider minimizing the control input (acceleration/deceleration though the gas/brake pedal position) that results in minimizing transient engine operation. If we minimize transient engine operation, we have direct benefits in fuel consumption \cite{Malikopoulos2008b} and emissions since internal combustion engines are optimized over steady state operating points (constant torque and speed)  \cite{Malikopoulos2010a,Malikopoulos2015a}. Fuel consumption 
was quantified by using the polynomial metamodel proposed in \cite{kamal_ecological_2011} which yields vehicle fuel consumption
as a function of speed, $v(t)$, and control input, $u(t)$, namely 
\begin{center}
\begin{equation}
\dot{f_{v}}=\dot{f}_{cruise}+\dot{f}_{accel},\label{eq:40}
\end{equation}
\par\end{center}
where $t\in\mathbb{R}^{+}$ is the time, $\dot{f}_{cruise}=w_{0}+w_{1}\cdot v(t)+w_{2}\cdot v^{2}(t)+w_{3}\cdot v^{3}(t)$
estimates the fuel consumed by a vehicle traveling at a constant speed
$v(t)$, and $\dot{f}_{accel}=u(t)\cdot(n_{0}+n_{1}\cdot v(t)+n_{2}\cdot v(t)^{2})$
is the additional fuel consumption caused by acceleration $u(t)$.
The polynomial coefficients $w_{n}$, $n=0,\,\ldots,\,3$ and $n_{s}$,
$s=0,\,1,\,2$ were calculated from experimental data. For the case
studies we considered in this paper, all vehicles were the same with
the parameters reported in \cite{kamal_ecological_2011},
where the vehicle mass was $M_{v}=1,200\,kg$, the drag coefficient
was $C_{D}=0.32$, the air density was $\rho_{a}=1.184\,km/m^{2}$,
the frontal area was $A_{F}=2.5\,m^{2}$, and the rolling resistance
coefficient was $\mu=0.015$.

\subsection{Testbed network}
\label{sec:4a} 
\begin{figure}[ptb]
\centering
\includegraphics[width=3 in]{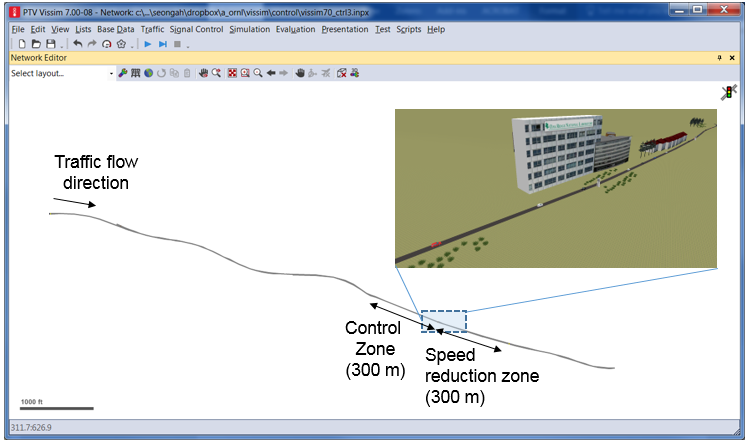}

\caption{Testbed network developed in VISSIM}

\label{fig:3}

\end{figure}

The testbed network consists of a single-lane corridor $2,000$ $m$ length long (Fig. \ref{fig:3}) that includes a speed reduction zone, $300$ $m$ length long, located 
downstream. As the capacity at the entrance of the speed reduction zone drops, the stop-and-go congestion is meant to generate under the baseline scenario even without feeding additional volumes through ramps.  
The speed limit inside the speed reduction zone is set to be $15.6$ $m/s$. We also designated a $300$ $m$ long control zone, right before the speed reduction zone. The length of the control zone was determined appropriately to provide adequate distance to address traffic congestion that would have occurred at the speed reduction zone (i.e., bottleneck) downstream. It is noted that neither the speed limit at the speed reduction zone, nor the length of the control zone affect the analytical, closed-form solution for system optimality. However, by varying the length of the control zone, or the speed limits, might yield different quantitative results on fuel consumption, travel time, and throughput.

It is noted that the VISSIM model was calibrated with a reference to the guideline of the Highway Capacity Manual 2010 \cite{HCM_2010}. According to this manual, the capacity of two-lane highways under the baseline scenario is indicated as $1,700$ $veh/h$. In this study, however, the design capacity was relaxed to $1,800$ $veh/h$ considering that the capacity increases in one-way versus a two-way highway as it has been shown based on empirical data \cite{ACTCanada2012,Spahr}. The minimum safe distance $\delta(t)$ in VISSIM which is defined
as the distance a driver would maintain while following vehicle
can be expressed as follows \cite{ptv_ptv_2014}

\begin{equation}
\delta(t)=c_{0}+c_{1}\cdot v_{ave}(t)\label{eq:41}
\end{equation}
where $c_{0}$ is the standstill distance between two vehicles,
$c_{1}$ is the headway time, and $v_{ave}(t)$ is the average speed.
In our case, we used the default value of VISSIM for $c_{0}$, e.g., $1.5$ $m,$ while $c_{1}$ was adjusted to $1.2$ $sec$,
thereby the maximum traffic flow was approximated at near the desired value of $1,800$ $veh/h$. In addition, some of the key parameters affecting VISSIM's car following behavior used include: the standstill distance of $1.5$ $m,$ headway time of 1.2 $sec$, and the absolute space difference between minimum and maximum gap distance as 4 $m$.

\subsection{Experimental set-up}

\label{sec:4b}To evaluate the effectiveness of the efficiency of the proposed approach
under varying traffic volume conditions, we considered the following cases: (i) traffic volume of $1,620$ $veh/h$ which is 10\% less
than the capacity, (ii) traffic volume of $1,800$ $veh/h$ at the capacity,
and (iii) traffic volume of $1,980$ $veh/h$ which is 10\% more than the
capacity. For all these cases, the total simulation time was
$1,000$ $sec$. Five replications of each simulation
case were conducted to account for the effect of stochastic components
of traffic and drivers' behaviors. It is noted that five replications produced statistically
significant results at a 95\% confidence level. 

The parameters used for the proposed control algorithm are summarized
in Table \ref{table:1}. According to a guideline published by the
Federal Highway Administration \cite{dowling_traffic_2004},
suggested maximum acceleration and deceleration are   
$3.04$ $m/s^{2}$ and  $-4.5$ $m/s^{2}$, respectively. Several recent studies have employed some relaxed values considering the advancement of the vehicle technologies, in which the thresholds ranged from $4.0$ $m/s^2$ to $8.0$ $m/s^2$ \cite{Lee:2013aa, Zhou:2012aa}. Since the vehicle considered here are automated, the
maximum acceleration threshold was assigned to be $4.5$ $m/s^{2}$
and the maximum deceleration was adopted as the guideline taking account
the safety and comfortable driving behavior. However, these values could be modified as necessary without any implications in the proposed framework.

\begin{table}

\caption{Control parameters}

\begin{centering}
\begin{tabular}{cc}
\hline 
Parameter & Value\tabularnewline
\hline 
Min. speed & 10 $m/s$\tabularnewline
Max. speed & 35$m/s$\tabularnewline
Max. acceleration & 4.5$m/s^{2}$\tabularnewline
Max. deceleration & -4.5$m/s^{2}$\tabularnewline
\hline 
\end{tabular}
\par\end{centering}
\label{table:1}

\end{table}

We considered the following three cases for comparison: (i) a baseline scenario associated with human drivers based on a car-following model, (ii) the state-of-the-art VSL algorithm, SPEed ControllIng ALgorIthm using Shock wave Theory (SPECIALIST) and (iii) the vehicular-based SPD-HARM algorithm proposed by the US DOT \cite{Ma:2016aa}. 

The baseline scenario was designed to emulate the human drivers' behavior using the Wiedemann model in VISSIM \cite{Wiedemann1974}. The comparison analysis with the baseline represents the net benefit of the control algorithm over the status that would have been resulted without the control algorithm. The SPECIALIST is a
proactive VSL algorithm that projects traffic conditions for a short-term horizon using MPC. It utilizes  shock wave theory
to generate control speed and duration, and
thus, it does not require complicated computation. In addition, it includes only a few
parameters with physical interpretations for feasible field
implementations. In this study, SPECIALIST was modeled
using C\# programming and implemented in the VISSIM using its
COM interface. Since SPECIALIST is based on a mesoscopic
model that utilizes the spot-based measurement collected at a fixed
location and aggregated for a certain period of time, detector stations
were evenly embedded at every $75$ $m$ along the corridor to estimate
the local traffic states. The traffic state of each detector station
was estimated every $60$ $sec$ by using the aggregated estimation
of the latest $60$ $sec$ interval, and the activation of VSL was examined
every $60$ $sec$ as well. 
The vehicles within the control zone
were ensured to follow the VSL control speed at 100\% compliance
rate without perception-reaction time. Such ideal condition was necessary
for a fair comparison with the proposed control algorithm which assumed
100\% automated vehicle penetration. The SPECIALIST algorithm had
several parameters that can be selected by the operator. For the best
performance of the algorithm, the parameters were tuned with several
iterations. The thresholds of  maximum speed and capacity were chosen
as $15.6$ $m/s$ and $1,500$ $veh/h$, respectively, which were determined after
empirical trials to find the minimum values where traffic congestion
was not observed under the VSL implemented at the 100\% automated
vehicle market penetration. 

A vehicular-based SPD-HARM algorithm, called Simple Speed Harmonization (SH) Algorithm \cite{Ma:2016aa}, was developed through the support of US DOT in an effort to realize field implementation of the SPD-HARM algorithm using CAVs. The simple SPD-HARM algorithm was created in inspiration from the work in \cite{Lu2014}; however, it was simplified to demonstrate the feasibility under the context of probe-based traffic management using CAVs, namely 
\begin{equation}
s_{i}(x,t)=m(t)\cdot x_{i}(t)+b(t) \label{eq:42}
\end{equation}
\begin{equation}
m(t)=\frac{s_{m}(t)-s_{n}(t)}{\Delta x_{mn}}\label{eq:43}
\end{equation}
\begin{equation}
b(t)=s_{m}(t)\label{eq:44}
\end{equation}
where $s_{i}(x,t)$ denotes the speed of vehicle $i$ at time $t$ and $x_{i}$ denotes the location of vehicle $i$ at time $t$. With successful validation of field implementation \cite{Ma:2016aa}, the effectiveness of the Simple SH algorithm in mitigating the traffic oscillations has been updated by Learn \cite{Learn2016}. The recent results showed that the probe-based SPD-HARM algorithm effectively mitigated the traffic oscillations compared to the base case of no-control, while the travel time and fuel consumption were increased over the base case \cite{Learn2016}.

\subsection{Results and Analysis}
\label{sec:4c}

\begin{figure}
\begin{center}
\includegraphics[width=3.5 in]{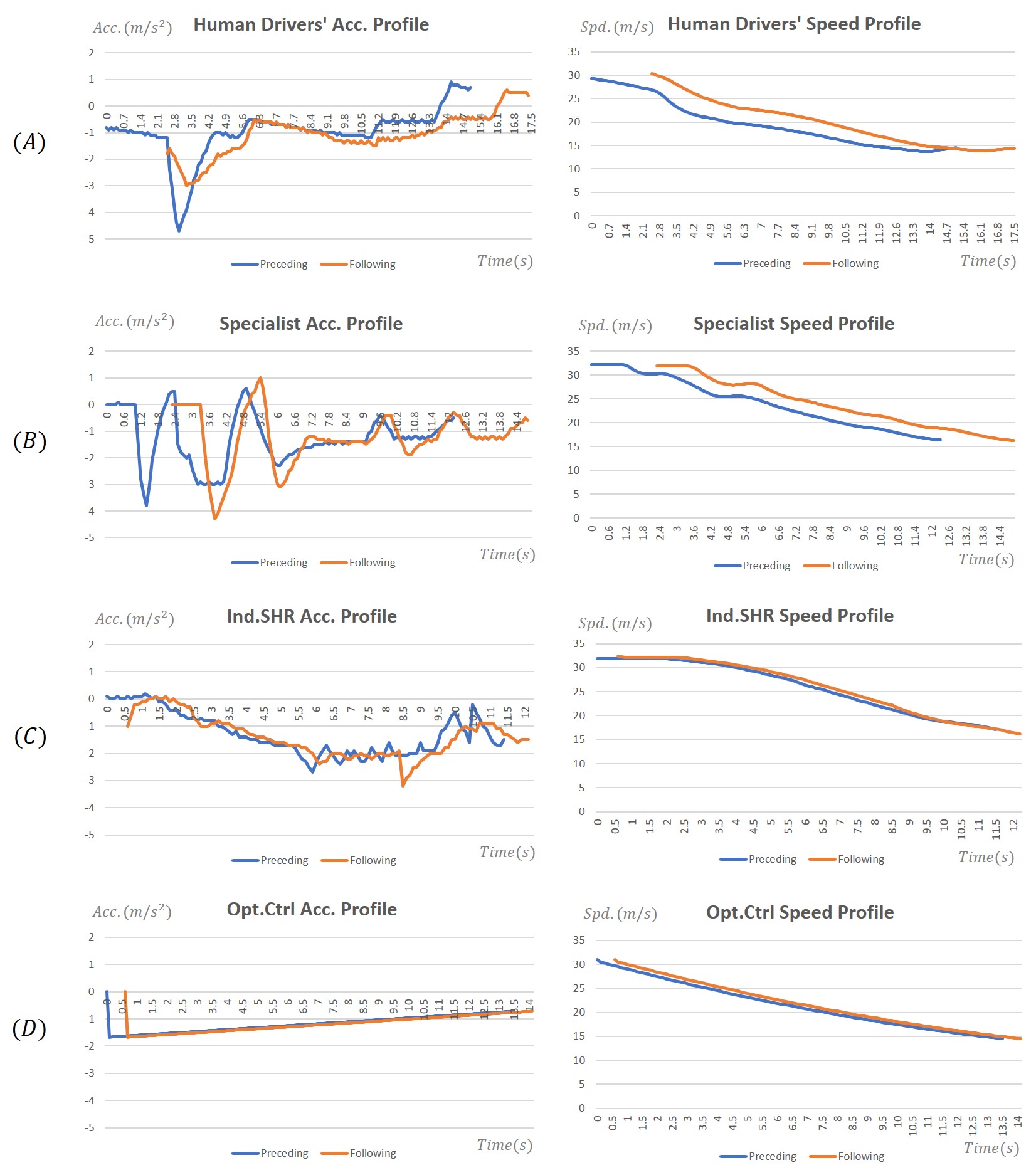}
\caption{Accelerations and speeds profiles of all cases}
\label{fig:4}
\end{center}
\end{figure}

We investigated acceleration and speed profiles of the four cases including the proposed control algorithm, VSL, vehicle-based SPD-HARM algorithm and human driven vehicles. Figure \ref{fig:4} shows the accelerations and speeds of two  adjacent vehicles traveling within the control zone. These trajectories were selected from the simulation testbed at the same time for all cases. 

While the speed trajectories of all cases show vehicles enter into the control zone at 31 m/s and exits the control zone around 15.6 m/s, the acceleration profiles are very different. As expected, human drivers case and the variable speed limit (VSL) case showed large fluctuations. The individual speed harmonization (SPD-HARM) case showed better acceleration profile without explicitly considering optimal control within the control zone. Finally, the proposed optimal control case demonstrates the most smooth acceleration profile, which is the linear function as designated by \eqref{eq:20}. 

\begin{figure}
\begin{center}
\includegraphics[width=3.5 in]{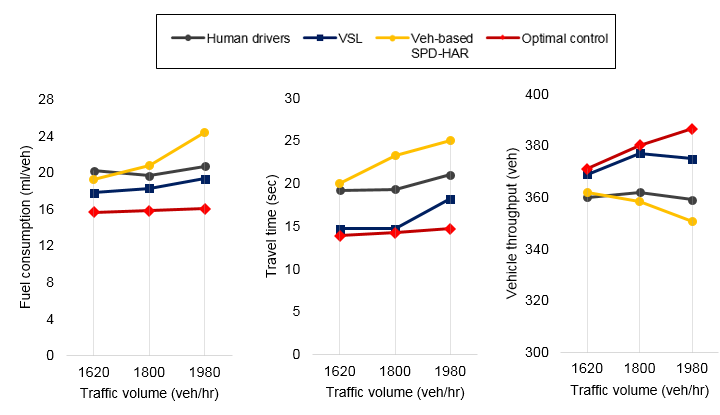}
\caption{Comparisons of the baseline scenario (human-driven vehicles), VSL, and optimal control algorithm.}
\label{fig:5}
\end{center}
\end{figure}

The results corresponding to fuel consumption, travel time, and throughput of the three approaches (i.e., VSL, vehicular-based SPD-HARM algorithm, and the proposed control algorithm) are shown in Fig. \ref{fig:5} and compared to a baseline scenario using human-driven vehicles. Fuel consumption per vehicle for the three approaches demonstrate similar patterns with the results of travel time as shown in Fig. \ref{fig:5} (a) and (b). By optimizing the vehicles' acceleration/deceleration inside the control zone, the time of reaching the speed reduction zone is controlled optimally, and thus the recovery time from the congested area is minimized, and therefore we have improvement in travel time as shown in Fig. \ref{fig:5} (b). Furthermore, each vehicle avoids getting into a stop-and-go driving mode, thereby conserving momentum and energy. Eliminating the vehicles' stop-and-go driving aims at minimizing  transient engine operation, and thus we have direct benefits in fuel consumption, as shown in Fig. \ref{fig:5} (a), since internal combustion engines are optimized over steady state operating points (constant torque and speed) \cite{Malikopoulos2008b}. 

The proposed control algorithm significantly reduces fuel consumption of each vehicle by 19-22\% over the baseline scenario, by 12-17\% over the VSL algorithm, and by 18-34\% over the vehicular-based SPD-HARM algorithm for the three traffic volume cases considered here. Using the vehicular-based SPD-HARM,  fuel consumption was statistically not different from that of the baseline scenario for the traffic volume cases  corresponding to the road capacity and the one of 10\% less than capacity. When the traffic volume was 10\% more than the capacity,  fuel consumption increased by 18\% over the baseline scenario. These observations are consistent with the results in an experiment conducted by a research group of FHWA \cite{Learn2016}. 

As shown in Fig. \ref{fig:5} -(b) and (c), the proposed control algorithm improved travel time and throughput for all traffic volumes over the baseline scenario, VSL, and the vehicular-based SPD-HARM algorithm.  In particular, travel time was improved by 26-30\% over the baseline scenario, by 3-19\% over the VSL algorithm, and by 31-39\% over the vehicular-based SPD-HARM algorithm for the three traffic volume cases. Both VSL and the proposed control algorithm reduced the travel time and improved the vehicle throughput under all three traffic volume cases. The proposed control algorithm reduced travel time and throughput by 19\% and 3\% respectively, compare to VSL under the traffic volume higher than the capacity. On the contrary, the vehicular-based SPD-HARM increased travel compared to the baseline scenario. This is a reasonable outcome since the vehicular-based SPD-HARM is designed to lead harmonizing speed throughout the traffic stream, but not optimizing its performance. However, the vehicle throughput at the bottleneck area was not significantly different than those of baseline scenario when the traffic volume was less than the capacity or at the capacity at 95\% confidence level as shown in Fig. \ref{fig:5}-(c). This implies that even though the vehicular-based SPD-HARM might increase the average travel time along the control zone, it eventually resulted in resolving potential shockwaves at the downstream bottleneck area and helped release the vehicles effectively unless the traffic congestion is initiated. 

To ensure statistical significance of the results, the t-tests were conducted for all three measurements by paring up of two scenarios out of the four (i.e., human drivers, VSL, vehicular-based SPD-HARM, and the optimal control algorithm). The statistical results showed that the p-values were close to zeros (rounded at the second decimal points), which implied that the measurement of certain scenario was significantly different from the other two scenarios at the 99th percentile confidence level. The simulation results were further assessed by identifying the minimum sample size using the Sample Size Determination Tool (Version 2.0) which was developed based on the FHWA sample size determination methodology \cite{VDOT:2013aa}. In the tool, the minimum sample size is determined as a function of the statistics of the initial simulation runs, confidence level and the tolerance error. In this paper, the 95th percentile confidence interval and 5\% error tolerance value were used as suggested in the manual \cite{VDOT:2013aa}. The initial five runs of simulation satisfy the minimum sample size for all scenarios per the tolerance level of 5\% with the 95th percentile confidence level.

It is important to highlight that the proposed control algorithm does not require a phase to clear congestion or a technique to prevent the bottleneck formation in order to improve vehicle throughput. Instead, the vehicle throughput was improved by having all individual vehicles proactively determines their optimal trajectories to the target location while the minimum spacing from the preceding vehicle is assured.


\section{Concluding Remarks and Future Research}
In the ``new world" of massive amounts of information from vehicles and infrastructure, what we used to model as uncertainty (noise or disturbance) becomes additional input or extra state information in a much higher-dimensional vector. The processing of such multiscale information requires new approaches in order to overcome the curse of dimensionality. Then the question becomes ``how much information do we need and what we can do with it?" It seems clear that the availability of this information has the potential to ease congestion, reduce energy usage, and  diminish traffic accidents by enabling vehicles to rapidly account for changes in their mutual environment. The approach we proposed in this paper demonstrated that automation and control can be used to improve traffic flow.

In particular, we addressed the problem of optimally controlling the speed of a number
of automated vehicles before they enter a speed reduction
zone on a freeway.  We formulated
the control problem and used Hamiltonian analysis to provide an analytical,
closed-form solution that can be implemented in real time. The
solution yields the optimal acceleration/deceleration of each vehicle
under the hard safety constraint of rear-end collision avoidance. The effectiveness of the proposed solution was demonstrated through
simulation and it was shown that the proposed approach can reduce 
both fuel consumption and travel time.

In our proposed framework, we considered that all vehicles are automated. Also we did not consider lane changing. An important direction for future research is to relax these assumptions and investigate the implications in the solution.
The assumption of perfect information might also impose barriers
in a potential implementation and deployment of the proposed framework.
Although it is relatively straightforward to extend our results in
the case that this assumption is relaxed, future research should investigate
the implications of having information with errors and/or delays to
the system behavior. Finally, considering lane changing and mixed
traffic (e.g., automated vehicles and human-driven vehicles) would
eventually aim at addressing the remaining practical consequences
of implementing this framework.


\section{ACKNOWLEDGMENTS}

This research was supported in part by ARPAE's NEXTCAR program under the award number DE-AR0000796 and in part by the SMART Mobility Initiative of the Department of Energy.  This research
project was also partially supported by the Global Research Laboratory
Program through the National Research Foundation of Korea (NRF) funded
by the Ministry of Science, ICT \& Future Planning (2013K1A1A2A02078326).
These supports are gratefully acknowledged.

\bibliographystyle{IEEETran}
\bibliography{TRB}

\begin{IEEEbiography}[{\includegraphics[width=1.1in,height=1.25in,clip,keepaspectratio]{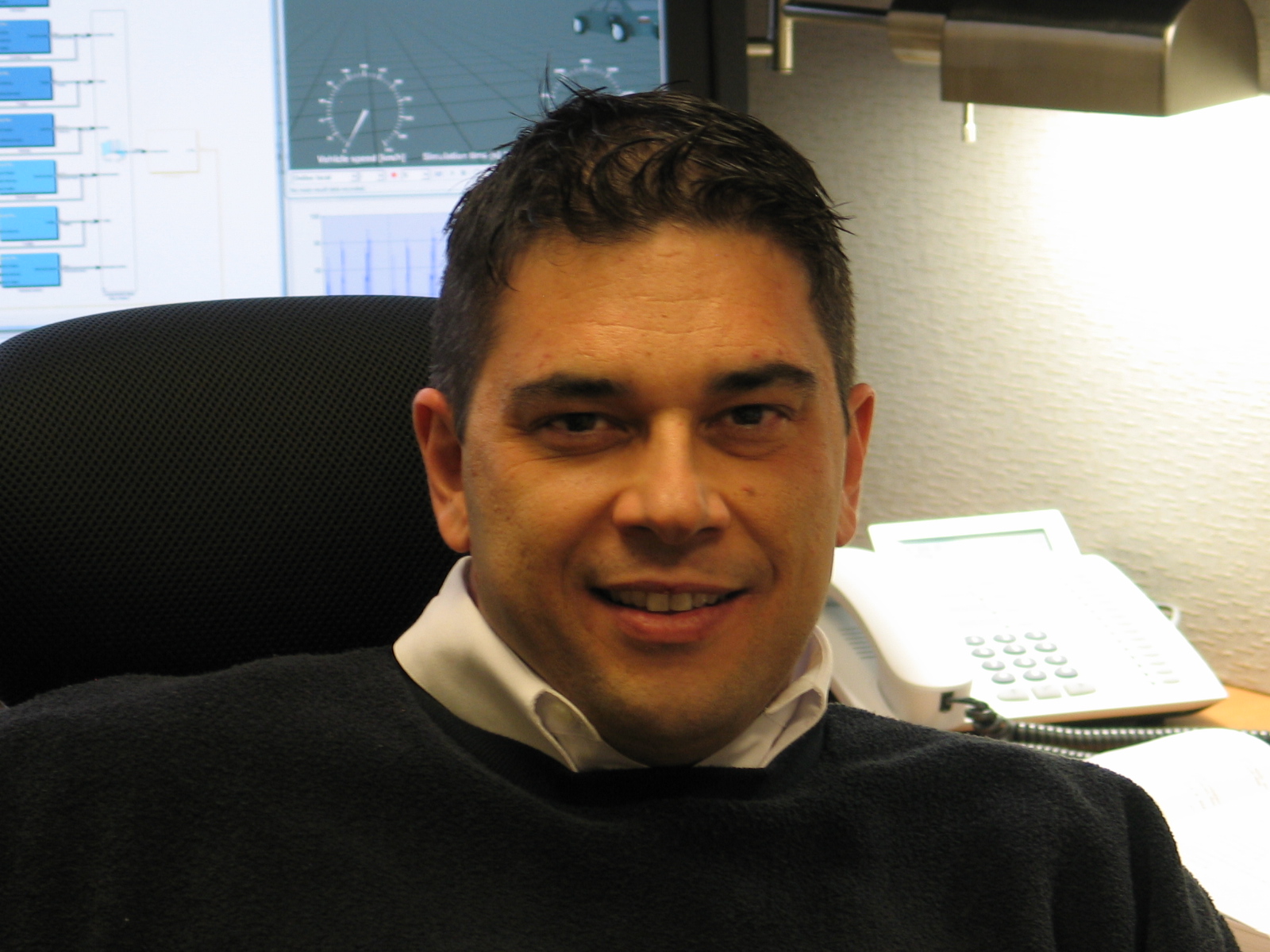}}]{Andreas A. Malikopoulos}
(M2006, SM2017) received a Diploma in Mechanical Engineering from the National Technical University of Athens, Greece, in 2000. He received M.S. and Ph.D. degrees from the Department of Mechanical Engineering at the University of Michigan, Ann Arbor, Michigan, USA, in 2004 and 2008, respectively. 
He is an Associate Professor in the Department of Mechanical Engineering at the University of Delaware (UD). Before he joined UD, he was the Deputy Director and the Lead of the Sustainable Mobility Theme of the Urban Dynamics Institute at Oak Ridge National Laboratory, and a Senior Researcher with General Motors Global Research \& Development. His research spans several fields, including analysis, optimization, and control of cyber-physical systems; decentralized systems; and stochastic scheduling and resource allocation problems. The emphasis is on applications related to sociotechnical systems, energy efficient mobility systems, and sustainable systems. He is currently an Associate Editor of the IEEE Transactions on Intelligent Vehicles and IEEE Transactions on Intelligent Transportation Systems. He is a member of SIAM and a Fellow of the ASME.
\end{IEEEbiography}

\begin{IEEEbiography}[{\includegraphics[width=1.1in,height=1.25in,clip,keepaspectratio]{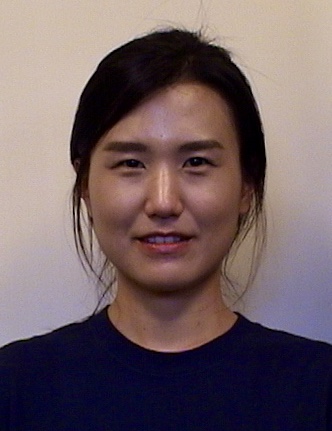}}]{Seongah Hong}
received the B.S. degree in Urban Planning and Engineering from Yonsei University, Seoul, Korea in 2012 and the M.S. from the Civil Engineering in the University of Minnesota, Duluth, MN in 2014, and a Ph.D. degree in the Civil and Environmental Engineering Department at the University of Virginia, Charlottesville, VA, USA. She has participated in various research project about traffic operations, intelligent transportation systems, connected automated vehicles applications, cyber-security and monitoring system, and sustainable traffic control devices. She has conducted extensive work of modeling and analyzing the effectiveness of the traffic operations applications with the customized control algorithms under diverse scenarios.
\end{IEEEbiography}

\begin{IEEEbiography}[{\includegraphics[width=1.1in,height=1.25in,clip,keepaspectratio]{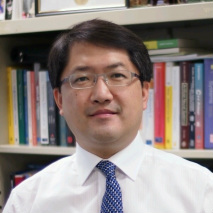}}]{B. Brian Park}
received the B.S. and the M.S. from the Hanyang University, Seoul, Korea, in 1993 and 1995, respectively, and the Ph.D. from the Texas A\&M University in 1998. Brian Park is an Associate Professor of Civil and Environmental Engineering Department at the University of Virginia. Prior to joining the University of Virginia, he was a Research Fellow at the National Institute of Statistical Sciences and a Post-Doctoral Research Associate at North Carolina State University. Dr. Park is a recipient of PTV America Best Paper Award, Outstanding Reviewer Award from the American Society of Civil Engineers, Jack H. Dillard Outstanding Paper Award from the Virginia Transportation Research Council and Charley V. Wootan Award (for best Ph.D. dissertation) from the Council of University Transportation Centers.  He is an Editor in Chief of the International Journal of Transportation, an Associate Editor of the American Society of Civil Engineers Journal of Transportation Engineering, Journal of Intelligent Transportation Systems and the KSCE Journal of Civil Engineering, and an editorial board member of the International Journal of Sustainable Transportation. Furthermore, he is a member of TRB (a division of the National Academies) Vehicle Highway Automation Committee and Artificial Intelligence and Advanced Computing Applications Committee, and chair of Simulation subcommittee of Traffic Signal Systems Committee. He is also Chair of Advanced Technologies Committee of ASCE Transportation and Development Institute.
\end{IEEEbiography}

\begin{IEEEbiography}[{\includegraphics[width=1.1in,height=1.25in,clip,keepaspectratio]{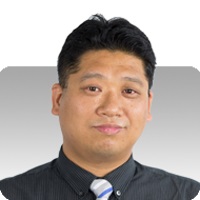}}]{Joyoung Lee}
received the B.S. degree in transportation engineering from Hanyang University, Ansan, Korea, in 2000 and the M.S. and Ph.D. degrees in civil and environmental engineering from the University of Virginia (UVA), Charlottesville, 2007 and 2010, respectively. He is currently an Assistant Professor with the John A. Reif, Jr. Department of Civil and Environmental Engineering at the New Jersey Institute of Technology (NJIT). Before joining NJIT in 2013, he was a laboratory manager of the Saxton Transportation Operations Laboratory (STOL) at the Federal Highway Administration (FHWA) Turner-Fairbank Highway Research Center. Dr. Lee is an Associate Editor of the KSCE Journal of Civil Engineering. He is a member of Transportation Research Board (TRB) Travel time, Speed, and Reliability (TTSR) Subcommittee.
\end{IEEEbiography}

\begin{IEEEbiography}[{\includegraphics[width=1.1in,height=1.25in,clip,keepaspectratio]{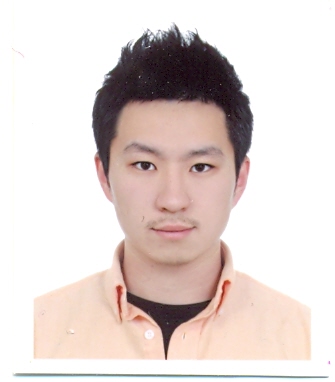}}]{Seunghan Ryu}
received the B.S. and M.S. degree in Transportation Engineering from Hanyang University, Ansan, Korea in 2011 and 2015, respectively. He is currently in his second year of Ph.D. program in the University of Virginia (UVA), Charlottesville, VA, USA. Before joining UVA, he was a researcher of The Korea Transportation Institute (KOTI), Sejong, Korea. He participated in various research project about traffic demand forecasting, demand modeling, feasibility study, intelligent transportation systems, connected automated vehicles applications, and vehicle platooning. 
\end{IEEEbiography}

\end{document}